\documentclass[10pt,twocolumn,twoside]{IEEEtran}
\usepackage{amsmath}
\usepackage{amssymb}            
\usepackage[pdftex]{graphicx}          
\usepackage[dvips]{epsfig}    
\newtheorem{coro}{Corollary}
\newtheorem{rem}{Remark}

\newtheorem{thm}{Theorem}
\newtheorem{prop}{Proposition}
\newtheorem{exmp}{Example}
\begin{document}
\title{Zero controllability in discrete-time structured systems}
\author{Jacob van der Woude
\thanks{Jacob van der Woude 
is with DIAM, EWI, Delft University of Technology, Mekelweg 4, 2628 CD, Delft, the Netherlands, (email: j.w.vanderwoude@tudelft.nl)}}
\maketitle
\begin{abstract}
In this paper we consider complex dynamical networks modeled by means of state space systems running in discrete time.
We assume that the dependency structure of the variables within the (nonlinear) network equations is known and use directed graphs to represent this structure.
The dependency structure also appears in the equations of a linearization of the network.
In order for such a linearization to be a good approximation of the original network, its state should stay as close as possible to the point of linearization.
In this paper, we investigate how the latter can be  achieved by an appropriate selection of states as driver nodes, so that through these driver nodes the whole state of the network can be steered to the point of linearization.
We present conditions in graph terms for this to be possible and derive an algorithm for the associated driver node selection, possibly of smallest size.

\end{abstract}
\newtheorem{algorithm}{Algorithm}
\section{Introduction and motivation}\label{introduction}
In this paper we consider complex dynamical networks modeled by means of state space systems running in discrete time.
We assume that we know how the equations describing the network look like, meaning that we know how the variables depend on each other.
More specifically, we assume that we have an explicit difference equation for each of the variables and knowledge on which variables play a role in each of these equations.
Hence, the dependency structure of the variables in the equations is assumed to be known.

In general, the network models are too large, too complicated and have too many unknown parameters, etc., to be studied or even to be used in practice.
For this reason, we assume that we can approximate the network equations by linear ones, for instance, by linearizing them.
In such linearizations, the equation structure follows from the dependency structure of variables within the equations of the original network.
However, the coefficients in these linearized equations in general will have an unknown value.
This is due to modeling errors, parameter uncertainties in the original network equations, etc.
Nevertheless, although not known, the value of these coefficients will in general be nonzero.
Hence, the result of the linearization of the network is a linear discrete-time model, with a known zero-nonzero coefficient structure for the equations, but with nonzero coefficients having an unknown value.

The obtained linear model is only valid in a neighborhood of the point of linearization.
However, in general the size of this neighborhood is unknown.
The linear model is a perfect approximation of the network when the state of the network is at the point of linearization or, said differently, when the state of the linearization is equal to zero.
Hence, the linear approximation is a good representation of the original network model when its state is close to the zero state.

To that end, it would be useful if we could steer the state of the linearization to the zero state, and preferably in finite time.
This is especially important in cases when this state is deviating from the zero state too much and may reach the (unknown) boundary of the region of linearization.
This steering can be achieved by controlling (some of) the states in the network.
Such states are then attached to controls and are called driver nodes in literature.
A big question then is which states should actually be selected to become a driver state.
Another question concerns the minimal number of driver nodes in order to have that the state can be steered to zero.

In this paper we will address the above questions and provide answers to them.
To that end, we recall some aspects of zero controllability for discrete-time systems and connect them to the above questions.

The outline of this paper is as follows.
In Section \ref{ProbFormul} we recall linear discrete-time systems and known results on zero controllability.
In Section \ref{systemsandgraphs} we introduce the type of structured systems that we consider in this paper together with a graph representation, and some notions and results from graph theory.
The introduced concepts are illustrated by means of an example.
In Section \ref{genericresults} we recall some of the well-known results on structural controllability for linear discrete-time structured systems.
In Section \ref{mainresults} we present our main results in terms of necessary and sufficient conditions for a linear discrete-time  structured system to be finite time zero controllable.
The conditions will be stated in graph theory terms.
In Section \ref{algorithm} we use our main result to obtain a method for driver node selection such that the state can be steered to zero in finite time.
Also the method will be illustrated by means of an example.
The paper is concluded by means of Section \ref{conclusions} containing remarks and topics for future research.

The current paper is partly based on the conference paper \cite{vdWoude:1997} where zero controllability (or deadbeat controllability) already was studied.
The driver node selection aspect in the context of zero controllability is however completely new.

\section{Problem formulation} \label{ProbFormul}
In order to address the questions raised in the introduction, we first consider the following class of linear discrete time systems
\begin{equation}\label{linsys}
x(k+1)=Ax(k)+Bu(k), \quad k=0,1,\ldots
\end{equation}
with state $x(k) \in \mathbb{R}^n$ and control $u(k) \in \mathbb{R}^m$.
The variable $k$ stands for the discrete time, and for the time being $A$ and $B$ are real numerically specified matrices of dimensions $n \times n$ and $n \times m$, respectively.

We recall that system (\ref{linsys}) is said to be controllable if the state of the system (\ref{linsys}) can be steered from any initial value $x_0$ at time $k=0$ to any final value value $x_1$ if there exists some finite time $k=\tau$ and a sequence of controls $u(k)$ on times $k=0,1,\ldots,\tau-1$ such that $x(\tau)=x_1$.

As a special case we call system (\ref{linsys}) is finite time zero controllable if the state of the system (\ref{linsys}) can be steered from any initial value $x_0$ at time $k=0$ to the zero state at some finite time $k=\tau$.

Later in this paper we assume that only the zero-nonzero structure of matrix $A$ is known and that matrix $B$ is not present.
The goal then is to find a matrix $B$ of certain structure, corresponding to a driver node selection, such that the zero controllability problem mentioned in the introduction can be solved.

It is easily seen from (\ref{linsys}) that
\begin{equation}\label{solution}
x(l)-A^l x(0) = \sum \limits_{k=0}^{l-1} A^{l-1-k} B u(k),
\end{equation}
for all integer $l \geq 0$.

The next list of well-known equivalent statements can be proved.

\begin{thm}\label{control}
Given system (\ref{linsys}) the following statements are equivalent.
\begin{enumerate}
\item[(i)] System (\ref{linsys}) is controllable,
\item[(ii)] $\mathbb{R}^n = \mbox{\rm im}\,\left(B, AB, \ldots, A^{n-1}B \right)$,
\item[(iii)] \mbox{\rm rank}\,$\left( A - \lambda I, B \right) = n$ for every eigenvalue $\lambda$ of matrix $A$,
\item[(iv)] \mbox{\rm rank}\,$\left( A - z I, B \right) = n$ for every $z \in \mathbb{C}$.
\end{enumerate}
\end{thm}

Using the expression in (\ref{solution}) with $l=n$ and the Cayley Hamilton theorem, the equivalence of statements $(i)$ and $(ii)$ can be easily proved.
The equivalence of statements $(iii)$ and $(iv)$ is immediate.
The equivalence of statements $(i)$ and $(iii)$ (or $(iv)$) is known as the Hautus test for controllablity.

In a similar way it can be proved that the statements in the following theorem are equivalent.

\begin{thm}\label{zerocontrol}
Given system (\ref{linsys}) the following statements are equivalent.
\begin{enumerate}
\item[(i)] System (\ref{linsys}) is zero controllable,
\item[(ii)] $\mbox{\rm im}\,{A}^n \subseteq \mbox{\rm im}\,\left(B, AB, \ldots, A^{n-1}B \right)$,
\item[(iii)] \mbox{\rm rank}\,$\left( A - \lambda I, B \right) = n$ for every nonzero eigenvalue $\lambda$ of matrix $A$,
\item[(iv)] \mbox{\rm rank}\,$\left( A - z I, B \right) = n$ for every $z \in \mathbb{C}, z \not =0$.
\end{enumerate}
\end{thm}

In fact, the equivalence of statements $(i)$ and $(ii)$ in Theorem \ref{zerocontrol} again follows easily from (\ref{solution}) with $l=n$ and the Cayley Hamilton theorem.
The equivalence of statements $(i)$ and $(iii)$ (or $(iv)$) follows from the Hautus test for stabilizability with only the origin in the complex plane as stability region.

From statement $(ii)$ in Theorem \ref{zerocontrol} it is clear that a sufficient condition for system (\ref{linsys}) to be zero controllable is that matrix $A$ is nilpotent, i.e., $A^n=0$.
In that case matrix $B$ can even be the zero matrix.

The latter is precisely the situation we mentioned in the introduction.
Indeed, later in this paper we assume that we only know the zero-nonzero structure of matrix $A$ and that matrix $B$ is not present.
Then we first check whether matrix $A$ is nilpotent in structural sense, meaning that this property only depends on the structure in $A$, and not on the numerical values within this matrix.
If it happens that matrix $A$ is nilpotent in structural sense, the system is generically zero controllable, even without using a control.ƒ

However, if matrix $A$ is not nilpotent in structural sense, we will investigate how we should choose driver nodes, defining (the structure of) matrix $B$, so that the system becomes generically zero controllable.
Also we will be able to determine the minimal number of driver nodes needed for generic zero controllability.

\section{Structured systems and directed graphs}\label{systemsandgraphs}
As indicated in the introduction, in many cases only the zero-nonzero structure of the matrices $A$ and $B$ is known.
The actual value of the nonzero entries in the matrices is then not known exactly.

The zero-nonzero structure of the system matrices $A$ and $B$ can be represented by means of a directed graph.
This directed graph will be denoted by ${\cal G}=({\cal V},{\cal E})$  with a {\em vertex set} ${\cal V}={\cal X} \cup {\cal U}$ and an {\em edge set} ${\cal E} = \{(x_j,x_i) | a_{ij} \not =0\} \cup \{(u_j,x_i) | b_{ij} \not =0\}$.
Here ${\cal X}$ denotes the set of {\em state vertices}, ${\cal U}$ the set of {\em input vertices}, $(x_j,x_i)$ denotes an edge from $x_j$ to $x_i$, $(u_j,x_i)$ denotes an edge from $u_j$ to $x_i$, $a_{ij} \not =0$ means that the $(i,j)$-th entry of $A$ is a nonzero, and similarly for $b_{ij} \not =0$.

Given two vertices $s, t \in {\cal V}$ we say that there is a {\em path} in ${\cal G}$ from $s$ to $t$, or that $t$ can be reached from $s$, if there are vertices $v_{i_0}, v_{i_1},\dots,v_{i_k} \in {\cal V}$ such that $v_{i_0}=s$, $v_{i_k}=t$ and $(v_{i_j},v_{i_{j+1}}) \in {\cal E}$ for $j=0,1,\ldots, k-1$.
The {\em begin vertex} of the path is $s$ and the {\em end vertex} is $t$.
The path then has {\em length} $k$.

We say that a state vertex is {\em reachable} from the input vertices ${\cal U}$, if there is a path from some input vertex to that state vertex.
Correspondingly, we call a state vertex {\em unreachable} from the input vertices ${\cal U}$, if the state vertex can not be reached from any input vertex.
A path that starts in ${\cal U}$ is also referred to as a ${\cal U}$-rooted path.

We also use the notion of {\em cycle}.
A cycle is a path of which only the begin and end vertex coincide.
Clearly, a cycle must have all its vertices in ${\cal X}$.
Hence, a cycle is a collection of edges $(x_{i_j},x_{i_{j+1}}) \in {\cal E}$ for  $j=0,1,\ldots, k-1$ such that  vertices $x_{i_0}, x_{i_1},\dots, x_{i_{k-1}} \in {\cal X}$ are mutually distinct and $x_{i_0}=x_{i_k}$.

Occasionally, $a_{ij}$ will be referred to as the weight of edge $(x_j,x_i)$, and similarly for $b_{ij}$.
The weight of a path is given by the product of the weights of the edges it consists of.
It can be show that the $(i,j)$-th entry of $A^k$, denoted $[A^k]_{ij}$, can be obtained as the sum of the weights of all paths from vertex $x_j$ to vertex $x_i$ of length $k$, and similarly for $[A^kB]_{ij}$, being the $(i,j)$-th entry of $A^kB$.

In this paper we also will consider the graph ${\cal G}$ without taking into account the set of input vertices ${\cal U}$.
Then, ${\cal V}={\cal X}$ and ${\cal E}=\{(x_j,x_i)| a_{ij} \not = 0\}$.

Having introduced the reachability between vertices in ${\cal G}$, we define a maximal strongly connected component to be a largest set of vertices ${\cal C}$ in ${\cal X}$, such that for any two vertices $s, t$ in ${\cal C}$ there is a path from $s$ to $t$, and a path from $t$ to $s$, where the paths may have any length, including zero.

We say that a maximal strongly connected component is nontrival if it contains at least one edge between its vertices, otherwise it is called trivial.
It follows easily that a trivial maximal strongly connected component must consist of just one vertex with no edges, i.e., no self loops.

The graph ${\cal G}$ can be partitioned into a collection of maximal strongly connected components that can be ordered as follows.
Given maximal strongly connected components ${\cal C}_1, {\cal C}_2$, we define/denote ${\cal C}_1 \prec {\cal C}_2$ if there is a path in ${\cal G}$ from a (=any) vertex in ${\cal C}_1$ to a (=any) vertex in ${\cal C}_2$.
The partitioning into maximal strongly connected components can be done by means of well-known and efficient algorithms.

\begin{exmp}\label{example1}
Consider the pair of structured matrices $(A,B)$ given by
\[
A=
\left( \begin{array}{ccccc}
\ast&\ast&\ast&0&0\\
\ast&0&0&\ast&0\\
0&0&0&\ast&\ast\\
0&0&0&0&0\\
0&0&0&0&\ast
\end{array} \right),
\quad
B=
\left( \begin{array}{c}
0\\
0\\
0\\
\ast\\
0
\end{array} \right),
\]
where the $\ast$'s denote unknown nonzero entries.
The graph of the corresponding system is displayed in Figure \ref{Fig0}.
\begin{figure}[htb]
 \centering
 \includegraphics[height=6cm,width=8cm]{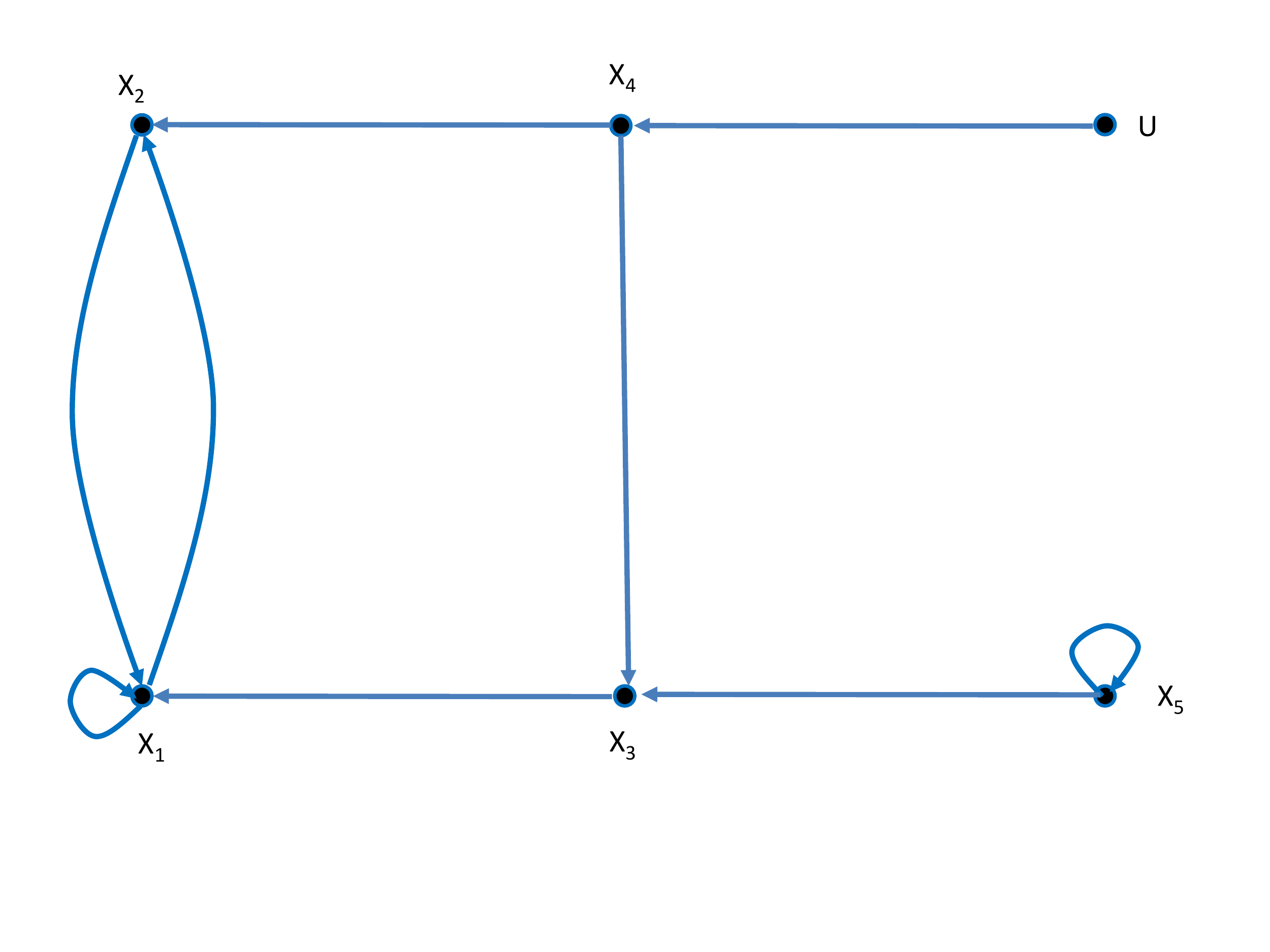}
 \vspace*{-1cm}
 \caption{Graph of example \ref{example1} \label{Fig0}}
\end{figure}

The collection of edges $\{(u,x_4),$ $(x_4,x_3),$ $(x_3,x_1),$ $(x_1,x_2)\}$ forms a ${\cal U}$-rooted path that ends in node $x_2$ and starts in ${\cal U}=\{u\}$.
The path has length $4$.
The weight of the path will occasionally be indicated by $a_{21}\, a_{13} \, a_{34} \, b_4$.
The path $(x_5,x_3),(x_3,x_1)$ from node $x_5$ to node $x_1$ has length $2$.

The collection of edges $\{(x_1,x_2),(x_2,x_1)\}$ forms a cycle in the state node set ${\cal X}=\{x_1,x_2,x_3,x_4,x_5\}$.
Also $\{(x_1,x_1)\}$ and $\{(x_5,x_5)\}$ are cycles in ${\cal X}$.

The ${\cal U}$-rooted path $\{(u,x_4),(x_4,x_3)\}$ and the cycles $\{(x_1,x_2),(x_2,x_1)\}$ and $\{(x_5,x_5)\}$ form a disjoint collection that covers every node of ${\cal X}$.

Note that the nodes $x_1, x_2, x_3$ and $x_4$ are reachable from ${\cal U}$. Node $x_5$ is not reachable from ${\cal U}$.

Ignoring the input node $u$, the maximal strongly connected components are ${\cal C}_1=\{x_1, x_2\}$, ${\cal C}_2=\{x_3\}$, ${\cal C}_3=\{x_4\}$ and ${\cal C}_4=\{x_5\}$.
The sets ${\cal C}_1$ and ${\cal C}_4$ are nontrivial ones, whereas ${\cal C}_2$ and ${\cal C}_3$ are trivial maximal strongly connected components.
The components are ordered as follows ${\cal C}_4 \prec {\cal C}_2$, ${\cal C}_3 \prec {\cal C}_2$ and ${\cal C}_2 \prec {\cal C}_1$, implying that also ${\cal C}_4 \prec {\cal C}_1$ and ${\cal C}_3 \prec {\cal C}_1$.
Note that the set $\{x_1\}$ is a strongly connected component, that however is not maximal, i.e., is not as large as possible.

Finally, note that entry $(1,5)$ of $A^3$ equals $a_{11} a_{13} a_{35}+a_{13} a_{35} a_{55}$.
\end{exmp}

\section{Generic results}\label{genericresults}
In this section we recall some results on generic controllability.
Also we present a result on the nilpotency of a square structured matrix.
Therefore, we assume that from now on that we only know the zero-nonzero structure of the matrices $A$ and $B$, forming a structured pair $\left(A,B\right)$ and defining a structured system of type (\ref{linsys}).

We say that a pair $(\bar{A},\bar{B})$ of numerically specified matrices is {\em admissible} to the structured pair $(A,B)$, if $(\bar{A},\bar{B})$ can be obtained from $(A,B)$ by fixing each of the nonzeros in the matrices $A$ and $B$ to a numerical value.

\subsection{Controllability}
We say that the structured pair $(A,B)$ is {\em generically controllable} if almost all admissible pairs $(\bar{A},\bar{B})$ are controllable in the usual sense, i.e., im $\left(\bar{B}, \bar{A}\bar{B},\ldots,\bar{A}^{n-1}\bar{B}\right)=\mathbb{R}^n$ or, equivalently, rank $\left( \bar{A} - \lambda I,\bar{B} \right) = n$ for every eigenvalue $\lambda$ of $\bar{A}$.
Here 'almost all' means 'all except for those in some proper algebraic variety (a set of zero measure) in the set of all admissible pairs'.
It can be shown that if there is one admissible pair that is controllable, then almost all admissible pairs are controllable, see \cite{lin:74}.

To present graph theoretic conditions for the generic controllability of a structured system we need two notions.

First, we say that the structured pair $(A,B)$ is {\em reducible}, if there is a permutation matrix $P$ such that
\[
P A P^{-1}=
\left(
\begin{array}{cc}A_{11}&A_{12}\\0&A_{22} \end{array}
\right), \quad
P B=
\left(
\begin{array}{c}B_{1}\\0 \end{array}
\right),
\]
with $A_{ij}$ an $n_i \times n_j$ matrix, $(i,j)=(1,1),(1,2),(2,2)$, $B_1$ an $n_1 \times m$ matrix and $n_1 \geq 0, n_2>0, n_1+n_2=n$.
If such a permutation does not exist, we say that the structured pair $(A,B)$ is {\em irreducible}.

Next, we define the {\em generic rank} of the structured $n \times (n+m)$ matrix $[A,B]$ to be the maximum of the ranks of the $n\times (n+m)$ matrices $[\bar{A},\bar{B}]$, taking into account all pairs $(\bar{A},\bar{B})$ that are admissible to the structured pair $(A,B)$.
We denote this rank by g-rank $[A,B]$.

We can now formulate the next well-known result, see \cite{DCV:03}, \cite{HM:79}, \cite{Reinschke:88}.

\begin{thm}\label{struccontr1}
The structured pair $(A,B)$, with $A \in \mathbb{R}^{n \times n}, B \in \mathbb{R}^{n \times m}$, is generically controllable if and only if the following two conditions hold.
\begin{enumerate}
\item[(i)] The structured pair $(A,B)$ is irreducible,
\item[(ii)] g-rank $[A,B]=n$.
\end{enumerate}
\end{thm}

Note that if the structured pair $(A,B)$ is reducible, the associated system (\ref{linsys}) can be rewritten as
\begin{eqnarray*}
x_1(k+1) &=& A_{11} x_1(k) + A_{12} x_2(k) + B_1 u(k),\\
x_2(k+1) &=& \hspace*{2.1cm} A_{22} x_2(k),
\end{eqnarray*}
where $x_1(k) \in \mathbb{R}^{n_1}$ and $x_2(k) \in \mathbb{R}^{n_2}$, with $n_1 \geq 0, n_2>0, n_1+n_2=n$, are subvectors of $x(k)$ obtained by a permutation of the components of $x(k)$.

Clearly, if the goal is to send each component of $x(k)$ to zero, this also has to apply to the components of $x_1(k)$ and $x_2(k)$.
For $x_2(k)$ this means that matrix $A_{22}$ should be such that this happens automatically for any initial condition.
Such behaviour will be investigated in more detail in the next subsection.

\subsection{Nilpotency}
We are going to focus on the original system (\ref{linsys}) without taking the control into account.
Hence, we are going to look at system
\begin{equation}\label{autsys}
x(k+1)=Ax(k),
\end{equation}
and investigate when the solution of this equation, i.e., $x(k)=A^k x(0)$, goes to zero in a finite number of steps.
Loosely, the latter means that $A$ has to be nilpotent or, equivalently, that $A$ only has the eigenvalue $0$.

Now recall that matrix $A$ is structured.
Following \cite{HM:79}, let $\nu(A)$ stand for the largest order of a principal minor of matrix $A$ that is not identically equal to zero.
According to Lemma 2 of \cite{HM:79}, matrix $A$ generically has $\nu(A)$ eigenvalues that are nonzero and mutually distinct, whereas the remaining $n-\nu(A)$ eigenvalues are all located at $0$.

Hence, it follows that $\nu(A)>0$ if and only if $A$ is generically not nilpotent.
Indeed, if $\nu(A)>0$, then $A$ generically has at least one eigenvalue that is nonzero, implying that $A$ is generically not nilpotent.
Conversely, if $A$ is generically not nilpotent, then $A$ generically has at least one nonzero eigenvalue, implying that $\nu(A)>0$.

From the proof of Lemma 2 in \cite{HM:79}, it follows that the characteristic polynomial of matrix $A$ has the form det$\,(sI_n-A)=s^{n-\nu(A)} p(s)$, where $p(s)$ is a polynomial with a generic degree $\nu(A)$ and mutually distinct zeros not located at zero.

The next result relates $\nu(A)=0$ to a property of the graph of matrix $A$.

\begin{thm}\label{genericnilpotent}
Let $A$ be a structured $n \times n$ matrix.
Then $\nu(A)=0$ if and only if the graph ${\cal G}$ of matrix $A$ contains no cycles.
\end{thm}

{\bf Proof}
$\Leftarrow$
Let $[A^k]_{ij}$ be the $(i,j)$-th entry of matrix $A^k$.
As indicated before, the entry is equal to the sum of the weights of all paths of length $k$ from vertex $x_j$ to vertex $x_i$.
Since ${\cal G}$ does not contain any cycle, there is a maximum to the length of any path in ${\cal G}$, say the maximum length is $l (\leq n)$.
It then follows that $[A^k]_{ij}$ with $k>l$ has to be zero.
Indeed, the entry can be seen as the sum of the weights of paths of length $k >l$ that each are zero because the paths can be thought of as being composed of edges of which at least one is nonexisting with zero weight.
Hence, it follows that $[A^k]_{ij}=0$ for all $k>l$ and all $1 \leq i,j \leq n$.
It even follows that the weight of the individual edges in ${\cal G}$ is irrelevant.
Hence, any matrix $\bar{A}$ that is admissible to $A$ is nilpotent and has only eigenvalues at $0$.
However, by Lemma 2 of \cite{HM:79}, the latter can only happen if $\nu(A)=0$.\\
$\Rightarrow$
Assume there exists a cycle in ${\cal G}$ of some positive length.
Now set the weight of the edges in the cycle to some (random) nonzero value and set the weight of all other edges equal to zero.
Then an admissible matrix $\bar{A}$ is found with a nonzero principal minor which easily can be shown to have a positive order, implying that $\nu(\bar{A})>0$.
The latter implies that $\nu(A)>0$, which contradicts the fact that $\nu(A)=0$.
$\hfill \Box$

%

\section{Generic zero controllability}\label{mainresults}

In this section we combine the results of the previous sections and we derive a method to determine whether or not a system is generically zero controllable.
First we define what we mean by generic zero controllability.

We say that the structured pair $(A,B)$ is generically zero controllable if almost all pairs $(\bar{A},\bar{B})$  admissible to $(A,B)$ are zero controllable in the usual sense, i.e., im $\bar{A}^n \subseteq \mbox{\rm im} \left(\bar{B}, \bar{A}\bar{B}, \ldots, \bar{A}^{n-1}\bar{B} \right)$ or, equivalently, rank $\left( \bar{A} - \lambda I,\bar{B} \right) = n$ for every nonzero eigenvalue $\lambda$ of $\bar{A}$.

Now consider a structured pair $(A,B)$ and assume that $(A,B)$ is not irreducible.
So, there exists a permutation matrix $P$ such that
\begin{equation}\label{decomp}
\left( P A P^{-1} \quad P B \right) =
\left(
\begin{array}{ccc}A_{11}&A_{12}&B_1\\0&A_{22}&0 \end{array}
\right),
\end{equation}
with $n_1 \geq 0, n_2>0, n_1+n_2=n$ and $(A_{11},B_1)$ irreducible.
The latter can always be accomplished.

It now clearly follows that g-rank $\left(A-zI_n,B\right)=n$ for all $z \not = 0$ if and only if
g-rank $\left(A_{11}-zI_{n_1},B_1\right)=n_1$
and
g-rank $\left(A_{22}-zI_{n_2}\right)=n_2$ for all $z \not = 0$.
The latter means that generically $A_{22}$ does not have any eigenvalues that are nonzero, which can only happen if and only if $\nu(A_{22})=0$.
Indeed, if $\nu(A_{22})>0$, then generically $A_{22}$ has at least one eigenvalue not equal to zero.
Conversely, if $A_{22}$ generically has a nonzero eigenvalue, then there is a principle minor of positive order, implying that $\nu(A_{22})>0$.

To come to our conditions we now recall one of the results of \cite{HM:79}.
See also \cite{MU:87} (proposition 14.5) and \cite{vdWoude:2000} (theorem 6.2), where extensions of the next result are discussed.

\begin{prop}
Let ${\cal G}$ be the graph of the structured system described by the pair $(A,B)$.
If $(A,B)$ is irreducible, then generically rank $\left( A - z I_n, B \right) =n$ for all $z \not = 0$.
\end{prop}

Recalling that in decomposition (\ref{decomp}) the pair $(A_{11},B_1)$ is irreducible, we therefore have with decomposition (\ref{decomp}) that g-rank $(A-zI_n,B)=n$ for all $z \not =0$ if and only if g-rank $(A_{22}-z I_{n_2})=n_2$ for all $z \not =0$.
The condition that g-rank $(A_{11}-zI_{n_1},B_1)=n_1$ for all $z \not =0$ is automatically fulfilled by the irreducibility of $(A_{11},B_1)$.
The condition that g-rank $(A_{22}-z I_{n_2})=n_2$ for all $z \not =0$ is equivalent to $\nu(A_{22})=0$, which by Theorem \ref{genericnilpotent} means that the part of the graph ${\cal G}$ that is not reachable from ${\cal U}$ does not contain any cycle.

To present results in graph terms we decompose ${\cal X}={\cal X}_r \cup {\cal X}_u$ with ${\cal X}_r \cap {\cal X}_u=\emptyset$, where ${\cal X}_r$ denotes the set of reachable states from ${\cal U}$ and ${\cal X}_u$ denotes the set of unreachable states from ${\cal U}$.
Next we remove from ${\cal G}$ all vertices in ${\cal X}_r$, together with all edges that have a vertex in ${\cal X}_r$ as begin or as end vertex.
We denote the graph that remains by ${\cal G}_u$.
Clearly ${\cal G}_u$ represents the part of ${\cal G}$ that is not reachable from ${\cal U}$.

We can summarize our observations as follows.

\begin{thm}\label{struczerocontr1}
Let ${\cal G}$ be the graph of the structured system described by the pair $(A,B)$ with ${\cal G}_u$ the part of ${\cal G}$ that is not reachable from ${\cal U}$.
The structured pair $(A,B)$ is generically zero controllable if and only if ${\cal G}_u$ contains no cycles.
\end{thm}

{\it Example 1. (cont.)}:
Note that in Example \ref{example1} the unreachable part ${\cal G}_u$ consists of ${\cal G}_u=\big(\{x_5\},$ $\{(x_5,x_5)\}\big)$.
Since it contains a cycle the system in Example \ref{example1} is not generically zero controllable.
However, when in Example \ref{example1} the entry $a_{55}$ is replaced by a fixed zero, then ${\cal G}_u$ does not contain a cycle anymore and the system does become generically zero controllable.

The above theorem is important in the development of an algorithm to select driver nodes to make a system become generically zero controllable.

\section{Driver node selection}\label{algorithm}
We consider a system of the form
\[
x(k+1)=Ax(k),
\]
and assume that matrix $A$ is a structured matrix.
The objective is to find a driver node selection such that with the induced structured matrix $B$ a structured system is obtained that is generically zero controllable.
To that end, we consider the index $\nu(A)$ as introduced before.

If $\nu(A)=0$, then the above objective is realized, because matrix $A$ is nilpotent in structural sense, and the associated structured system is generically zero controllable automatically, even without using control.
By Theorem \ref{genericnilpotent} it can be checked easily whether or not $\nu(A)=0$.

If $\nu(A)>0$, then by Theorem \ref{genericnilpotent}, the graph contains cycles and can be decomposed into maximal strongly connected components of which at least one is nontrivial, i.e., contains a cycle.
Then to achieve the above objective some driver nodes should be selected, defining matrix $B$, such that the associated system becomes generically zero controllable.

Recall that if the graph of matrix $A$ contains cycles, then they are contained in the nontrivial maximal strongly connected components.
Hence, by Theorem \ref{struczerocontr1}, a driver node selection to achieve generic zero controllability should have an unreachable part that does not contain any of the nontrivial components or, equivalently, should have a reachable part that contains all nontrivial components.
Trivial components may appear in the unreachable part, since they don't contain any cycle.

To formulate the latter in graph terms, let ${\cal G}$ be the graph of the structured system described by matrix $A$ with state vertex set ${\cal X}$.
Let ${\cal D}$ be a given driver node selection and decompose ${\cal X}$ as ${\cal X}_r^{\cal D} \cup {\cal X}_u^{\cal D}$ with ${\cal X}_r^{\cal D}$ the set of states that are reachable from ${\cal D}$ and ${\cal X}_u^{\cal D}$ the other vertices.
Like before, remove from ${\cal G}$ all nodes in ${\cal X}_r^{\cal D}$ and all edges with begin or end node in ${\cal X}_r^{\cal D}$.
Denote the the remaining graph by ${\cal G}_u^{\cal D}$.
It can be thought of as the part of ${\cal G}$ that is not reachable from ${\cal D}$.
Then the application of Theorem \ref{struczerocontr1} can formulated as follows.

\begin{coro}\label{corr1}
Consider the graph ${\cal G}$ of matrix $A$ and let ${\cal D}$ be a driver node selection.
Denote ${\cal G}_u^{\cal D}$ for the part of ${\cal G}$ that is not reachable from ${\cal D}$.
Then the system is generically zero controllable using the driver nodes in ${\cal D}$ if and only if ${\cal G}_u^{\cal D}$ contains no cycle.
\end{coro}

To achieve generic zero controllability a driver node set can be selected as follows.

\begin{algorithm}
Given the structured system with matrix $A$ and associated graph ${\cal G}$.
\begin{itemize}
\item If ${\cal G}$ contains no cycles, then the system is generically zero controllable automatically.
\item If ${\cal G}$ does contain cycles, search for a driver node selection ${\cal D}$ such that every cycle is reachable from one of the nodes in ${\cal D}$.
For a minimal driver node selection look for the smallest of such driver node selections.
\end{itemize}
\end{algorithm}

\begin{rem}
The minimal number of driver nodes to solve the problem of generic zero controllability equals the minimal number of nodes from which all nontrivial components in the associated graph can be reached.
Because of the aperiodic structure in which the (non)trivial maximal strongly connected components can be arranged the computation of a minimal driver node set can be done by means of efficient methods.
\end{rem}

\begin{exmp}\label{example2}
Consider the system $x(k+1)=Ax(k)$ with the structured matrix
\[
A=
\left( \begin{array}{ccccccccccc}
0&0&\ast&\ast&0&0&0&0&0&0&0\\
\ast&0&0&0&0&\ast&0&0&0&0&0\\
0&\ast&0&0&\ast&0&0&0&0&0&0\\
0&0&0&\ast&0&0&0&0&0&0&0\\
0&0&0&0&\ast&0&0&\ast&0&0&0\\
0&0&0&0&0&0&\ast&0&0&0&0\\
0&0&0&0&0&\ast&0&\ast&0&0&0\\
0&0&0&0&0&0&0&0&\ast&\ast&\ast\\
0&0&0&0&0&0&0&0&0&0&0\\
0&0&0&0&0&0&0&0&0&0&0\\
0&0&0&0&0&0&0&0&0&0&0
\end{array} \right),
\]
where the $\ast$'s denote unknown nonzero entries.
The graph of the corresponding system is displayed in Figure \ref{Fig1} in which also the trivial and nontrivial maximal strongly connected components are indicated.
\begin{figure}[htb]
 \centering
 \includegraphics[width=8cm]{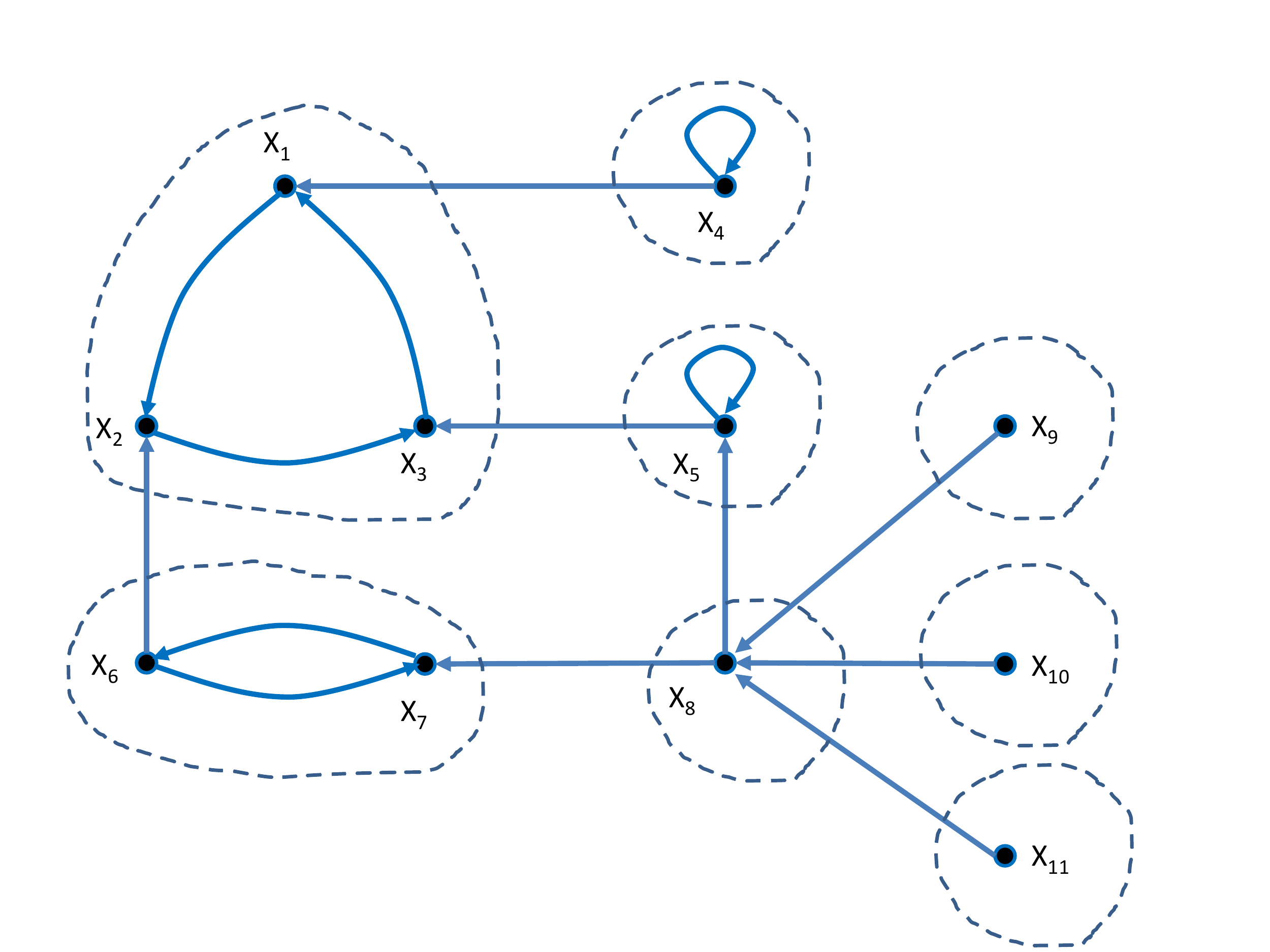}
 \vspace*{-1cm}
 \caption{Graph of Example 2 \label{Fig1}}
\end{figure}

Taking ${\cal D}=\{x_1, x_5\}$ as a driver node selection, it follows that the conditions in Corollary \ref{corr1} are not met.
Indeed, in the associated subgraph ${\cal G}_u^{\cal D}$ there is for instance a cycle at node $x_4$, i.e., a self loop $\{(x_4,x_4)\}$.
Furthermore, this unreachable graph contains the cycle $\{(x_6,x_7),(x_7,x_6)\}$.
Also ${\cal D}=\{x_1, x_5, x_7\}$ is not yet an adequate driver node selection, because the associated subgraph ${\cal G}_u^{\cal D}$ still contains a cycle at node $x_4$.

It is easily seen that any driver node set that meets the conditions of Corollary \ref{corr1} should contain node $x_4$.
Further, nodes $x_5, x_6, x_7$ should be at least reachable from a suitable set of driver nodes.
Note that through these nodes then also nodes $x_1, x_2, x_3$ are reachable from the selected driver nodes.
So, for instance, ${\cal D}=\{x_4, x_5, x_6\}$ is a suitable driver node set to make the system generically zero controllable, since all cycles can be reached from these nodes.

A smaller suitable driver node is ${\cal D}=\{x_4, x_8\}$, or alternatively ${\cal D}=\{x_4, x_k\}$ with $k=9,10,11$.
All make that the system becomes zero controllable, and all contain the minimal number of driver nodes to do so.
\end{exmp}

Occasionally the relation between driver nodes and controls is defined in different ways.
One definition requires that one control can be used to steer all driver nodes.
Another definition requires each driver node to have its own control.
This difference leads to different matrices $B$ to achieve generic zero controllability.
For the minimal driver nodes selection ${\cal D}=\{x_4,x_8\}$, given in Example \ref{example2}, the previous two definitions lead to the following versions of the structured matrix $B$
\[
B=
\left( \begin{array}{c}
0\\
0\\
0\\
\ast\\
0\\
0\\
0\\
\ast\\
0\\
0\\
0
\end{array} \right), \quad
B=
\left( \begin{array}{cc}
0&0\\
0&0\\
0&0\\
\ast&0\\
0&0\\
0&0\\
0&0\\
0&\ast\\
0&0\\
0&0\\
0&0
\end{array} \right),
\]
respectively.

 \section{Concluding remarks}\label{conclusions}
The idea underlying this paper is to keep a linear approximation of a complex real-life discrete time network valid as long as possible.
This can be done by keeping the state of the associated linearization at zero or having that it goes to zero in a finite number of steps.
The latter happens automatically for systems of which the matrix of the associated linearization is structurally nilpotent.
The graph of such systems does not contain any cycle.
For systems where the previous is not the case state nodes may have to be selected as driver nodes to control the entire system.
The objective then is to select these driver nodes such that the state of the linearization can be steered to zero in finite time.
The latter can be achieved if the driver nodes are selected in such a way that all cycles in the graph of the system can be reached from these nodes.
The selection of these nodes and the minimal number to achieve the desired behaviour is described in this paper.
From the paper it follows that for achieving generic zero controllability the connectability of the graph of the system is most relevant.
The latter is in contrast with recent results on controllability of complex networks where the focus is more on the generic rank of the structured system matrices, cf.\ \cite{Liu:11}.
A topic for future research may be a classification of nodes into essential, useful and useless driver nodes for achieving generic zero controllability.

\section{Acknowledgement}
The topic of this paper is inspired by the bachelor thesis of Tim den Haan, cf.\ \cite{TRBdH16},
who is greatly acknowledged for the nice collaboration during his project.
Further, the author wants to thank Christian Commault and Taha Boukhobza for useful remarks made in the process of writing this paper.
\bibliographystyle{plain}
\bibliography{Bib-ZeroContr}
\end{document}